\documentclass[10pt]{article}
\font\smallit=cmti10
\font\smalltt=cmtt10

\usepackage{amssymb,latexsym,amsmath,epsfig,amsthm,mathrsfs}
\usepackage{hyperref}
\usepackage{todonotes}
\usepackage{tabularray}
\UseTblrLibrary{diagbox}

\makeatletter

\renewcommand\section{\@startsection {section}{1}{\z@}
	{-30pt \@plus -1ex \@minus -.2ex}
	{2.3ex \@plus.2ex}
	{\normalfont\normalsize\bfseries\boldmath}}

\renewcommand\subsection{\@startsection{subsection}{2}{\z@}
	{-3.25ex\@plus -1ex \@minus -.2ex}
	{1.5ex \@plus .2ex}
	{\normalfont\normalsize\bfseries\boldmath}}

\renewcommand{\@seccntformat}[1]{\csname the#1\endcsname. }

\makeatother

\def\chopsticks/{\textsc{chopsticks}}
\def\octopus/{\textsc{simple chopsticks}}
\newcommand{\hands}[3]{(#1~|~#2)_{#3}}
\newcommand{\handstwo}[2]{(#1~|~#2)_{2}}
\newcommand{\outcomeN}{\mathcal{N}}
\newcommand{\outcomeP}{\mathcal{P}}
\newcommand{\outcomeL}{\mathcal{L}}
\newcommand{\outcomeR}{\mathcal{R}}

\newtheorem{theorem}{Theorem}
\newtheorem{lemma}{Lemma}
\newtheorem{proposition}{Proposition}

\newtheorem{notation}[theorem]{Notation}
\theoremstyle{definition}
\newtheorem{definition}{Definition}

\newtheorem{observation}[theorem]{Observation}

\emergencystretch=\maxdimen
\begin{document}
\begin{center}
	\uppercase{\bf {\sc simple chopsticks}: playing with any number of hands and fingers}
	\vskip 20pt
	{\bf Antoine Dailly}\\
	{\smallit Université Clermont-Auvergne, CNRS, Mines de Saint-Étienne, Clermont-Auvergne-INP, LIMOS, 63000 Clermont-Ferrand, France}\\
	{\tt antoine.dailly@uca.fr}\\
	\vskip 10pt
	{\bf Valentin Gledel}\\
	{\smallit Université Savoie Mont Blanc, CNRS UMR5127, LAMA, Chambéry, France}\\
	{\tt valentin.gledel@univ-smb.fr}\\
	\vskip 10pt
	{\bf Richard J. Nowakowski}\\
	{\smallit Dalhousie University, Canada}\\
	{\tt r.nowakowski@dal.ca}\\
	\vskip 10pt
	{\bf Carlos Pereira dos Santos}\\
	{\smallit Center for Mathematics and Applications (NovaMath), FCT NOVA, Portugal}\\
	{\tt cmf.santos@fct.unl.pt}\\
\end{center}
\vskip 20pt
\centerline{\smallit Received: , Revised: , Accepted: , Published: } 
\vskip 30pt

\centerline{\bf Abstract}
\noindent
\chopsticks/ is a game played by two players where they start with one finger raised on each hand. On their turn, each player moves by pointing an attacking hand at one of their opponent's hands. The number of fingers on the pointed hand increases by the number of fingers on the attacking hand. If, after a move, a hand contains more than five fingers, it is removed from play. There are also other rules that allow players to move fingers from one hand to another, but we focus on this simple setup.

We introduce a generalization of \chopsticks/, called \octopus/, in which the players may have any number of $n$-fingered hands. We find that having more hands than your opponent is generally good, and use this fact to fully characterize the outcomes of \octopus/ in the case where the players have 2-fingered hands.

\pagestyle{myheadings}
\markright{\smalltt INTEGERS: 24 (2023)\hfill}
\thispagestyle{empty}
\baselineskip=12.875pt
\vskip 30pt

\section{Introduction}\label{sec1}
\chopsticks/ is a well-known hand game \cite{Wikia2023,Wikib2023} with lots of different variants and rules. We are going to analyze the simplest version, where the players can only attack each other. Two players hold up their hands, originally with one finger raised in each. The players then take turns pointing one of their hands to one of their opponent's hands. The number of raised fingers in the pointed hand increases by the number of raised fingers in the attacking hand; if this number is more than its number of fingers\footnote{Usually five, for human beings.}, then the hand is removed from play. The first player to have all their hands\footnote{Again, usually two for human beings.} removed from play loses.

In this setup, the positions are combinatorial short games\footnote{Each position has finitely many distinct subpositions and admits no infinite run.}, which can be fully analyzed in the case of human beings. However, during the Combinatorial Game Theory Colloquium II \cite{CGTC2017}, in January 2017, a generalization of the game was proposed, so that beings with more (or less!) than two hands and five fingers could also analyze the winning strategies in their version of \chopsticks/.

We propose a generalization of \chopsticks/, called \octopus/, in which a position is fully expressed by two lists (describing the number of raised fingers in the hands of each player) and a \emph{finger count} (which is the highest possible number of fingers in the hands of the players). Note that, in this generalization, although the two players may have different numbers of hands, all the hands have the same number of fingers. The formal definition of \octopus/ is the following.

\begin{definition}
	Let $\ell$ and $r$ denote the number of Left's and Right's hands, respectively. Let $n$ be the finger count, \emph{i.e.}, the number of fingers on each hand of each player. Let $(x_1,\ldots,x_\ell)$ and $(y_1,\ldots,y_r)$ be two non-decreasing sequences of positive integers such that $x_i,y_i \leqslant n$ for all $i$.
	
	A position $\hands{x_1,\ldots,x_\ell}{y_1,\ldots,y_r}{n}$ of \octopus/ has the following options:
	\begin{itemize}
		\item The options for Left are the positions $\hands{x_1,\ldots,x_\ell}{y_1,\ldots, y_j+x_i ,\ldots,y_r}{n}$ for every $i \in \{1,\ldots,\ell\}$ and $j \in \{1,\ldots,r\}$ (if $y_j+x_i>n$, then, it is removed from the list);
		
		\item The options for Right are the positions $\hands{x_1,\ldots,x_i+y_j ,\ldots,x_\ell}{y_1,\ldots,y_r}{n}$ for every $i \in \{1,\ldots,\ell\}$ and $j \in \{1,\ldots,r\}$ (if $x_i+y_j>n$, then, it is removed from the list).
	\end{itemize}
\end{definition}

\begin{notation}
	Let $P=\hands{x_1,\ldots,x_{\ell}}{y_1,\ldots,y_r}{n}$ be a position of \octopus/. The Left-option $\hands{x_1,\ldots,x_{\ell}}{y_1,\ldots, y_j+x_i ,\ldots,y_r}{n}$ is denoted by $P(x_i \rightarrow y_j)$ or simply $x_i \rightarrow y_j$, if it is settled that the initial position is $P$. Similarly, the Right-option $\hands{x_1,\ldots,x_i+y_j ,\ldots,x_{\ell}}{y_1,\ldots,y_r}{n}$ is denoted by $P(y_j \rightarrow x_i)$ or $y_j \rightarrow x_i$.
\end{notation}

\begin{notation}
	If $k,a$ are two positive integers, $k^a$ represents the sequence where $k$ is repeated $a$ times.
\end{notation}

We assume that the reader is acquainted with the basic concepts of Combinatorial Game Theory (CGT) as presented in any of~\cite{AlberNW2007,BerleCG1982,Con1976,Siegel2013}. We only consider normal play where the player who cannot move loses.

In this paper, we study some configurations of \octopus/. First, in Section~\ref{sec2}, we look at the case where each player has only one hand, and find that, while the outcomes are trivially predictable, the game values exhibit a strange behaviour. We then give, in Subsection~\ref{subsec31}, a general result, which states that having more hands than the opponent is generally a good thing. Next, in Subsection~\ref{subsec32}, we fully characterize the outcomes of the \octopus/ positions where the players' hands have two fingers each. Finally, although a general theory for the case where the players' hands have three fingers is difficult to obtain, we show in Subsection~\ref{subsec33} that the initial position in that case is an $\outcomeN$-position.

\section{When each player has one hand}\label{sec2}
When each player has only one hand, computing the outcome is easy, but the game values are still interesting and non-trivial. Since each player only has one hand, they have only one option. Moreover, the initial positions are symmetric  and therefore the only outcomes are Next player win, $\outcomeN$, and Previous player win, $\outcomeP$.

\begin{theorem}\label{thm:11n} Let $f_1=1$, $f_2=1$, and $f_{i+2}=f_{i}+f_{i+1}$ be the \textit{Fibonacci} numbers. Let the positive integer $n$ be the finger count. If there exists an integer $i$ such that $f_{2i}\leqslant n<f_{2i+1}$ then $\hands{1}{1}{n}$ is a Next player win, otherwise it is a Previous player win.
\end{theorem}
\begin{proof} By rules, the starting position is $\hands{1}{1}{n}$. Hence, since $f_1=f_2=1$, the starting position is $\hands{f_2}{f_1}{n}$. If Left moves first, after her move, the position is $\hands{f_2}{f_1+f_2}{n}=\hands{f_2}{f_3}{n}$. And, by induction, after Left's $i$th move, the position is $\hands{f_{2i}}{f_{2i+1}}{n}$. Analogously, after Right's $i$th move,
	the position is $\hands{f_{2i+2}}{f_{2i+1}}{n}$. Thus, if $f_{2i}\leqslant n <f_{2i+1}$, then Left wins, and, if $f_{2i+1}\leqslant n <f_{2i+2}$, then Left loses.
	
	The result now follows since if Right moves first, then the analysis is the same.
\end{proof}

\noindent\textbf{Game values of  $\hands{i}{j}{n}$ and $\hands{1}{1}{n}$.}
Regarding game values, it is natural to assume that a move that removes the last opponent's hand from the game is a move to $0=\{\,|\,\}$. With that assumption, \octopus/ is a \emph{dicotic} ruleset, which means that, at each stage of play, either
both players or neither player can move~\cite{AlberNW2007,BerleCG1982,Con1976,Siegel2013}.

It only makes sense to determine game values if there are situations where the disjunctive sum appears naturally. That is the case for \octopus/. Suppose that children are playing in teams, for example, five kids against five kids. They can form five pairs, each consisting of two players from different teams. The rules are as follows: during their turn, a team chooses one pair and the member of their team of that pair makes a move. The goal is to make a move such that at the end of your team's turn, all of the opponent team's players with remaining hands have no opponent with remaining hands. With this setup, each pair is a disjoint component. Another interpretation can be made if we think of \octopus/ as a board game. In this interpretation, we have an $n\times n$ squared chessboard, and each singular piece is placed on a square. Left can move a piece upwards a number of squares equal to the column where it is positioned. Right can move a piece to the right a number of squares equal to the row where it is positioned. If the move causes the piece to go outside the board, that move is terminal as far as that piece is concerned. Thus, from $(a,b)$, Left can move a piece to $(a+b,b)$ and Right to $(a,a+b)$. If $a+b>n$, the piece is removed from the game. Naturally, this game can be played with $k$ non-conflicting pieces, with the first player to remove all pieces from the board winning. If $k>1$, then we have a disjunctive sum, and the game values are the key to understanding the positions.

The game values are surprisingly complicated. We are only able to report on facts obtained with help of the computer program CGSuite~\cite{CGSuite2023}. Table~\ref{tab-conwayChopN20} shows the game values for $\hands{i}{j}{20}$ with $1\leqslant i,j\leqslant 20$ (smaller finger counts exhibits similar behavior). 
It is interesting to observe that all values in the table are \emph{uptimals}~\cite{McKay2007} (they can be expressed as sums of $*$ and powers of $\uparrow$); even the value $J=\{0\,|\,\uparrow^{[2]}\}$, which mysteriously has special occurrence in \octopus/, is the uptimal $.22*=\uparrow+\uparrow+\uparrow^2+\uparrow^2+*$. One could conjecture that all game values of \octopus/ are uptimals. However, that is not true. For example, $\hands{2}{5}{40}=\{J\,|\,0\}$ and $\hands{1}{1}{12}=\pm J$ are not uptimals.

\begin{table}[htb]
	\scalebox{0.6}{
		\begin{tblr}{
				colspec={ccccccccccccccccccccc},
				rows={m},
				vlines,
				hlines,
				hline{1,21,22} = {2pt},
				vline{1,2,22} = {2pt}
			}
			
			$20$ & $*$ & $*$ & $*$ & $*$ & $*$ & $*$ & $*$ & $*$ & $*$ & $*$ & $*$ & $*$ & $*$ & $*$ & $*$ & $*$ & $*$ & $*$ & $*$ & $*$ \\
			
			$19$ & $0$ & $*$ & $*$ & $*$ & $*$ & $*$ & $*$ & $*$ & $*$ & $*$ & $*$ & $*$ & $*$ & $*$ & $*$ & $*$ & $*$ & $*$ & $*$ & $*$ \\
			
			$18$ & $\uparrow$ & $0$ & $*$ & $*$ & $*$ & $*$ & $*$ & $*$ & $*$ & $*$ & $*$ & $*$ & $*$ & $*$ & $*$ & $*$ & $*$ & $*$ & $*$ & $*$ \\
			
			$17$ & $\uparrow^{[2]}$ & $0$ & $0$ & $*$ & $*$ & $*$ & $*$ & $*$ & $*$ & $*$ & $*$ & $*$ & $*$ & $*$ & $*$ & $*$ & $*$ & $*$ & $*$ & $*$ \\
			
			$16$ & $\uparrow^{[3]}$ & $\uparrow$ & $0$ & $0$ & $*$ & $*$ & $*$ & $*$ & $*$ & $*$ & $*$ & $*$ & $*$ & $*$ & $*$ & $*$ & $*$ & $*$ & $*$ & $*$ \\
			
			$15$ & $\uparrow^{[4]}$ & $\uparrow$ & $0$ & $0$ & $0$ & $*$ & $*$ & $*$ & $*$ & $*$ & $*$ & $*$ & $*$ & $*$ & $*$ & $*$ & $*$ & $*$ & $*$ & $*$ \\
			
			$14$ & $\uparrow^{[5]}$ & $\uparrow^{[2]}$ & $\uparrow$ & $0$ & $0$ & $0$ & $*$ & $*$ & $*$ & $*$ & $*$ & $*$ & $*$ & $*$ & $*$ & $*$ & $*$ & $*$ & $*$ & $*$ \\
			
			$13$ & $\uparrow^{[6]}$ & $\uparrow^{[2]}$ & $\uparrow$ & $0$ & $0$ & $0$ & $0$ & $*$ & $*$ & $*$ & $*$ & $*$ & $*$ & $*$ & $*$ & $*$ & $*$ & $*$ & $*$ & $*$ \\
			
			$12$ & $\uparrow^{[7]}$ & $\uparrow^{[3]}$ & $\uparrow$ & $\uparrow$ & $0$ & $0$ & $0$ & $0$ & $*$ & $*$ & $*$ & $*$ & $*$ & $*$ & $*$ & $*$ & $*$ & $*$ & $*$ & $*$ \\
			
			$11$ & $\uparrow^{[8]}$ & $\uparrow^{[3]}$ & $\uparrow^{[2]}$ & $\uparrow$ & $0$ & $0$ & $0$ & $0$ & $0$ & $*$ & $*$ & $*$ & $*$ & $*$ & $*$ & $*$ & $*$ & $*$ & $*$ & $*$ \\
			
			$10$ & $\uparrow^{[9]}$ & $\uparrow^{[4]}$ & $\uparrow^{[2]}$ & $\uparrow$ & $\uparrow$ & $0$ & $0$ & $0$ & $0$ & $0$ & $*$ & $*$ & $*$ & $*$ & $*$ & $*$ & $*$ & $*$ & $*$ & $*$ \\
			
			$9$ & $*$ & $*$ & $\uparrow^{[2]}$ & $\uparrow$ & $\uparrow$ & $0$ & $0$ & $0$ & $0$ & $0$ & $0$ & $*$ & $*$ & $*$ & $*$ & $*$ & $*$ & $*$ & $*$ & $*$ \\
			
			$8$ & $\downarrow$ & $*$ & $*$ & $*$ & $\uparrow$ & $\uparrow$ & $0$ & $0$ & $0$ & $0$ & $0$ & $0$ & $*$ & $*$ & $*$ & $*$ & $*$ & $*$ & $*$ & $*$ \\
			
			$7$ & $\Downarrow*$ & $\downarrow$ & $*$ & $*$ & $*$ & $*$ & $0$ & $0$ & $0$ & $0$ & $0$ & $0$ & $0$ & $*$ & $*$ & $*$ & $*$ & $*$ & $*$ & $*$ \\
			
			$6$ & $0$ & $\downarrow_{[2]}$ & $*$ & $*$ & $*$ & $*$ & $*$ & $\downarrow$ & $0$ & $0$ & $0$ & $0$ & $0$ & $0$ & $*$ & $*$ & $*$ & $*$ & $*$ & $*$ \\
			
			$5$ & $\uparrow$ & $0$ & $\downarrow_{[2]}$ & $*$ & $*$ & $*$ & $*$ & $\downarrow$ & $\downarrow$ & $\downarrow$ & $0$ & $0$ & $0$ & $0$ & $0$ & $*$ & $*$ & $*$ & $*$ & $*$ \\
			
			$4$ & $\uparrow^{[2]}$ & $0$ & $0$ & $0$ & $*$ & $*$ & $*$ & $*$ & $\downarrow$ & $\downarrow$ & $\downarrow$ & $\downarrow$ & $0$ & $0$ & $0$ & $0$ & $*$ & $*$ & $*$ & $*$ \\
			
			$3$ & $*$ & $\{0\,|\,\uparrow^{[2]}\}$ & $0$ & $0$ & $\uparrow^{[2]}$ & $*$ & $*$ & $*$ & $\downarrow_{[2]}$ & $\downarrow_{[2]}$ & $\downarrow_{[2]}$ & $\downarrow$ & $\downarrow$ & $\downarrow$ & $0$ & $0$ & $0$ & $*$ & $*$ & $*$ \\
			
			$2$ & $\downarrow_{[2]}$ & $*$ & $\{\downarrow_{[2]}\,|\,0\}$ & $0$ & $0$ & $\uparrow^{[2]}$ & $\uparrow$ & $*$ & $*$ & $\downarrow_{[4]}$ & $\downarrow_{[3]}$ & $\downarrow_{[3]}$ & $\downarrow_{[2]}$ & $\downarrow_{[2]}$ & $\downarrow$ & $\downarrow$ & $0$ & $0$ & $*$ & $*$ \\
			
			$1$ & $0$ & $\uparrow^{[2]}$ & $*$ & $\downarrow_{[2]}$ & $\downarrow$ & $0$ & $\Uparrow*$ & $\uparrow$ & $*$ & $\downarrow_{[9]}$ & $\downarrow_{[8]}$ & $\downarrow_{[7]}$ & $\downarrow_{[6]}$ & $\downarrow_{[5]}$ & $\downarrow_{[4]}$ & $\downarrow_{[3]}$ & $\downarrow_{[2]}$ & $\downarrow$ & $0$ & $*$ \\
			
			\diagbox[dir=NE,linewidth=2pt]{$i$}{$j$}& $1$ & $2$ & $3$ & $4$ & $5$ & $6$ & $7$ & $8$ & $9$ & $10$ & $11$ & $12$ & $13$ & $14$ & $15$ & $16$ & $17$ & $18$ & $19$ & $20$ \\
			
	\end{tblr}}
	\caption{Game values for all positions $\hands{i}{j}{20}$ ($1\leqslant i,j\leqslant 20$).}
	\label{tab-conwayChopN20}
\end{table}

Table~\ref{tab-conwayChop11N} shows the  values and outcomes for $\hands{1}{1}{n}$ with $1\leqslant n\leqslant 606$. We notice several points:

\begin{enumerate}
	\item The positions can take three values: 0 for $\outcomeP$-positions and $*$ and $\pm J$ for $\outcomeN$-positions;
	\item The $\outcomeN$ and $\outcomeP$-positions switch when $n$ reaches the next fibonacci number (Theorem \ref{thm:11n});
	\item The value $\pm J$ appears when $n$ takes values in the sequence of the row sums of an unsigned characteristic triangle for the Fibonacci numbers (from 12 and beyond), described as the A110035 sequence in OEIS\footnote{\url{http://oeis.org/A110035}}, and until the next $\outcomeP$ position.
\end{enumerate}

We conjecture that this behaviour continues for the rest of the integers.

\begin{table}[htb]
	\centering
	\begin{tblr}{
			colspec={ccc},
			vlines,
			hline{1,2,20}={solid}
		}
		
		Interval of $n$ & Value of $\hands{1}{1}{n}$ & Outcome \\ 
		1 & $*$ & $\outcomeN$ \\
		2 & $0$ & $\outcomeP$ \\
		3--4 & $*$ & $\outcomeN$ \\
		5--7 & $0$ & $\outcomeP$ \\
		8--11 & $*$ & \SetCell[r=2]{m} $\outcomeN$ \\
		12 & $\pm J$ & \\
		13--20 & $0$ & $\outcomeP$ \\
		21--30 & $*$ & \SetCell[r=2]{m} $\outcomeN$ \\
		31--33 & $\pm J$ & \\
		34--54 & $0$ & $\outcomeP$ \\
		55--79 & $*$ & \SetCell[r=2]{m} $\outcomeN$ \\
		80--88 & $\pm J$ & \\
		89--143 & $0$ & $\outcomeP$ \\
		144--208 & $*$ & \SetCell[r=2]{m} $\outcomeN$ \\
		209--232 & $\pm J$ & \\
		233--376 & $0$ & $\outcomeP$ \\
		377--545 & $*$ & \SetCell[r=2]{m} $\outcomeN$ \\
		546--606 & $\pm J$ & \\ 
	\end{tblr}
	\caption{Game values of $\hands{1}{1}{n}$ for $n\leqslant 606$.}
	\label{tab-conwayChop11N}
\end{table}

\section{When both player have many hands} \label{sec3}

\subsection{A general result}\label{subsec31}

The following result states that having more hands than your opponent is generally good.

\begin{theorem}
	\label{thm-moreHands}
	Let $\ell$, $r$ and $n$ be positive integers, and $P=\hands{x_1,\ldots,x_{\ell}}{y_1,\ldots,y_r}{n}$ be a position of \octopus/. If $\ell > r$ (resp. $\ell<r$) and Left (resp. Right) is the first player, then Left (resp. Right) has a winning strategy.
\end{theorem}

\begin{proof}
	Assume that $\ell>r$, and let $P=\hands{x_1,\ldots,x_{\ell}}{y_1,\ldots,y_r}{n}$ be a position of \octopus/. By inducting on the options, we prove that Left has a winning strategy as the first player on $P$.  Note that since $\ell>r$, we have $\ell \geqslant 2$, \emph{i.e.}, Left has at least two options.
	
	Firstly, if Left has a move that removes an integer from Right's sequence (Right ``loses one hand''), then there exist $i$ and $j$ such that $x_i+y_j>n$. In all of Right's options of $P(x_i \rightarrow y_j)$, Left has more numbers in her list. Therefore, by induction, Left wins against all of Right's responses. Hence, $P(x_i \rightarrow y_j)$ is a winning move for Left.
	
	We may now assume that $x_i+y_j \leqslant n$ for all $i,j$, that is, Left cannot remove a Right's hand from the game.
	
	To begin with, suppose that Left has a ``safe move'' that prevents Right from answering with a removal of one of her hands, that is, a move  that prevents Right from removing any number from Left's sequence with his answer. In other words, let us assume that $x_1+y_1+x_\ell \leqslant n$. Again, against all Right replies in $P(x_1 \rightarrow y_1)$, Left maintains more numbers in her list.  Hence, by induction, Left wins and $P(x_1 \rightarrow y_1)$ is a winning move.

	Finally, assume that $x_1+y_1+x_\ell > n$, that is, if Left chooses $x_1 \rightarrow y_1$, Right has an answer that removes a Left's hand from the game. We claim that $x_\ell \rightarrow y_1$ is a winning move for Left. Note that Right must choose a response that eliminates one of Left's hands from play, otherwise an inductive argument, similar to those used in previous cases, shows that Left wins. Thus, suppose that Right removes an integer $x_i$ from Left's sequence. Left still has at least one term, say $x_j$, in her sequence. We know that $x_j \geqslant x_1$ and therefore $x_j+y_1+x_\ell > n$.  Left playing $x_j \rightarrow x_\ell+y_1$ removes $x_\ell+y_1$ from Right's sequence. Now, by assumption, for each $r,s$,  $x_r+y_s \leqslant n$, and against all Right moves, Left has more numbers in her list.  Hence, by induction, Left wins and $x_\ell \rightarrow y_1$ is a winning move for Left.
\end{proof}

\begin{lemma}
	\label{lem-lessHands}
	Let $P=\hands{x_1,\ldots,x_{\ell}}{y_1,\ldots,y_r}{n}$ be a \octopus/ position. Assume $\ell< r$ and Left is the first player. If Left has a winning strategy, then $\ell=r-1$ and the winning moves reduce the number of hands for Right.
\end{lemma}

\begin{proof}We cannot have $\ell<r-1$. If this were to happen, even if the first move was Left removing a number from Right's list, by Theorem \ref{thm-moreHands}, Right would have a winning response. This would contradict the fact that Left has a winning strategy. Hence, we have $\ell=r-1$. On the other hand, all winning moves must remove a number from Right's list, otherwise a similar contradiction is reached.\end{proof}

\subsection{Finger count 2}\label{subsec32}

We fully characterize the outcomes of the \octopus/ positions game with finger count~2. First, we use Theorem~\ref{thm-moreHands} to characterize the positions where one player has at least two more hands than their opponent.

\begin{proposition}
	\label{prop-magicNum2-wayMoreHands}
	Let $a,b,c,d$ be nonnegative integers such that $|(a+b)-(c+d)| \geqslant 2$. The \octopus/ position $\handstwo{1^a,2^b}{1^c,2^d}$ has outcome
	\begin{itemize}
		\item $\outcomeL$ if $a+b \geqslant c+d+2$;
		\item $\outcomeR$ if $c+d \geqslant a+b+2$.
	\end{itemize}
\end{proposition}

\begin{proof}
	This is a direct corollary of Theorem~\ref{thm-moreHands}.
\end{proof}

Thus, there are only two kinds of positions left to consider: those where both players have the same number of hands, and those where the first player has exactly one hand less than their opponent. We begin our study with the first case.

\begin{proposition}
	\label{prop-magicNum2-sameNumberOfHands}
	Let $a,b,c,d$ be nonnegative integers such that $a+b=c+d>0$. The \octopus/ position $\handstwo{1^a,2^b}{1^c,2^d}$ has outcome
	\begin{itemize}
		\item $\outcomeP$ if $b=d=0$, and $a$ and $c$ are odd;
		\item $\outcomeL$ if $b=0$, $d \neq 0$ and $c$ is odd;
		\item $\outcomeR$ if $d=0$, $b \neq 0$ and $a$ is odd;
		\item $\outcomeN$ in all the other cases.
	\end{itemize}
\end{proposition}

\begin{proof}
	For a \octopus/ position $P = \handstwo{1^a,2^b}{1^c,2^d}$, we define $f(P)$ as the quantity $2(a+c)+(b+d)$. This quantity will always decrease during the play, and serves as a type of ``position rank'' which can be used to establish an assumption of minimality.
	
	We prove the result by contradiction. Let $P = \handstwo{1^a,2^b}{1^c,2^d}$ (with $a+b=c+d$) be a position that does not verify the statement and that minimizes $f(P)$. First, note that we must have $a+b \geqslant 2$. Indeed, if $a+b=1$, then we have two cases: either $P \in \{\handstwo{1}{2},\handstwo{2}{1},\handstwo{2}{2}\}$ and it is an $\outcomeN$-position (the first player removes the hand of the other player) or $P = \handstwo{1}{1}$ and it is a $\outcomeP$-position (it has $\handstwo{1}{2}$ as a Left-option and $\handstwo{2}{1}$ as a Right-option). Both cases verify the statement.
	
	There are several possible cases to consider, depending on the possible values of $a$, $b$, $c$ and $d$. We are going to explore those different cases, showing that, in each of them, $P$ actually verifies the statement, leading to a contradiction.

	\begin{itemize}
		\item If $b=d=0$, then the first player (without loss of generality, we can assume that it is Left; the argument for Right is symmetric) will play to $P'=\handstwo{1^a}{1^{a-1},2}$. Note that $f(P')=2(2a-1)+1<2(2a)=f(P)$, and, by minimality of $f(P)$, the position $P'$ verifies the statement. If $a$ is even, then $P'$ is an $\outcomeL$-position, and thus Left wins on $P$ as the first player. Since, by symmetry, Right wins on $P$ as the first player, this implies that $P$ is an $\outcomeN$-position and verifies the statement; a contradiction. If $a$ is odd, then $P'$ is an $\outcomeN$-position, and hence $P$ is a $\outcomeP$-position and verifies the statement; again a contradiction.
		
		\item Assume now that $b=0$, $d \neq 0$, and $c$ is odd, so $P=\handstwo{1^a}{1^c,2^d}$. Observe that $a=c+d$, and, since $c$ is odd, $a>0$. Also, if $d-1=0$, then $a$ is even. First, suppose that Right is the first player. Right has two options: $P'=\handstwo{1^{a-1},2}{1^c,2^d}$ and $P'' = \handstwo{1^{a-1}}{1^c,2^d}$. On one hand, it is easy to check that $f(P')<f(P)$, thus $P'$ is an $\outcomeN$-position by minimality of $f(P)$; on the other hand, Left can play from $P''$ to $\handstwo{1^{a-1}}{1^c,2^{d-1}}$ which is either an $\outcomeL$- or a $\outcomeP$-position by minimality of $f(P)$. Thus, Left has a winning strategy if Right is the first player. Therefore, in order not to contradict the fact that $P$ does not satisfy the statement, playing first in $P$, Left must lose.
		
		However, if Left is the first player, then she can play to $P' = \handstwo{1^a}{1^c,2^{d-1}}$. If $d-1=0$, then we have $P' = \handstwo{1^a}{1^c}$ with $a>c$, and, after Right's reply, Left wins by Theorem \ref{thm-moreHands}. Otherwise, by Lemma~\ref{lem-lessHands}, the only possible winning move for Right is to $P'' = \handstwo{1^{a-1}}{1^c,2^{d-1}}$. But, by minimality of $f(P)$, $P''$ is an $\outcomeL$-position. In both cases Left wins, and we have a contradiction.
		
		The case where $d=0$, $b \neq 0$, and $a$ is odd is symmetric.
		
		\item Assume now that $a=0$. From $P$, the first player will necessarily remove one hand from their opponent, and the second
		player too, reaching a position $P''$ where $a$ is still equal to~0, and both players have one hand less than in $P$. By minimality of $f(P)$, since we cannot have $a=0$ and $b=0$, we are in the ``all the other cases'' of the statement. Hence, the position $P''$ is an $\outcomeN$-position, and thus this is also the case for $P$; a contradiction.
		
		The case where $c=0$ is symmetric.
		
		\item Assume now that $b=0$, $d \neq 0$, and $c$ is even and positive, so $P=\handstwo{1^a}{1^c,2^d}$. Observe that $a=c+d$ and, since $d \neq 0$, $a>0$. If Left is the first player, then she can play to $P' = \handstwo{1^a}{1^c,2^{d-1}}$. If $d-1=0$, then we have $P' = \handstwo{1^a}{1^c}$ with $a>c$, and, after Right's reply, Left wins by Theorem~\ref{thm-moreHands}. Otherwise, by Lemma~\ref{lem-lessHands}, we conclude that the only possible winning move for Right is to $P'' = \handstwo{1^{a-1}}{1^c,2^{d-1}}$. But, by minimality of $f(P)$, $P''$ is an $\outcomeN$-position, so Left has a winning strategy from $P$ as the first player. Furthermore, if Right is the first player, then he moves to $P'=\handstwo{1^{a-1}}{1^c,2^d}$. Once more, by Lemma~\ref{lem-lessHands}, the only possible winning move for Left is to $P''=\handstwo{1^{a-1}}{1^c,2^{d-1}}$. But, by minimality of of $f(P)$, $P''$ is an $\outcomeN$-position, so Right has a winning strategy from $P$ as the first player. Hence, $P$ is an $\outcomeN$-position; a contradiction.
		
		The case where $d=0$, $b \neq 0$, and $a$ is even and positive is symmetric.
		
		\item The only remaining possibility is that $a$, $b$, $c$ and $d$ are all positive integers in the position $P = \handstwo{1^a,2^b}{1^c,2^d}$. Assume without loss of generality that Left is the first player (if Right is the first player, then we can apply a symmetric argument). Left plays to $P'=\handstwo{1^a,2^b}{1^{c-1},2^d}$. By Lemma~\ref{lem-lessHands}, there are two possible winning moves for Right, which require to remove a Left's hand. However, in both cases, by minimality of $f(P)$, the position reached after Right's move will be either an $\outcomeL$-position (if $b=1$, $d \neq 0$, $c-1$ is odd, and Right played to $\handstwo{1^a}{1^{c-1},2^d}$) or an $\outcomeN$-position in all the other cases). Thus, Right will leave a winning position for Left, who was the first player on $P$. By symmetry, if Right is the first player on $P$, he also has a winning strategy, and thus $P$ is an $\outcomeN$-position; a contradiction.
	\end{itemize}
	
	The proof is complete.
\end{proof}

\begin{observation}
	Note that the $\outcomeP$-position corresponding to $a=b=c=d=0$ was left out of the statement of Proposition \ref{prop-magicNum2-sameNumberOfHands}. The reason for this is that it is a \emph{Garden of Eden}, \emph{i.e.}, a position that has no predecessor. There is no possible move that leads to this trivial position.
\end{observation}

We now study the positions where the number of hands of the players differ by exactly one.

\begin{proposition}
	\label{prop-magicNum2-oneMoreHand}
	Let $a,b,c,d$ be nonnegative integers such that $|(a+b)-(c+d)| = 1$.
	\begin{enumerate}
		\item If $a+b>c+d$, then the position $\handstwo{1^a,2^b}{1^c,2^d}$ of \octopus/ has outcome
		\begin{itemize}
			\item $\outcomeP$ if $c=d=0$;
			\item $\outcomeN$ if $d = 0$, $b > 0$ and $a$ is odd;
			\item $\outcomeL$ otherwise.
		\end{itemize}
		
		\item If $a+b<c+d$, then the position $\handstwo{1^a,2^b}{1^c,2^d}$ of \octopus/ has outcome
		\begin{itemize}
			\item $\outcomeP$ if $a=b=0$;
			\item $\outcomeN$ if $b = 0$, $d > 0$ and $c$ is odd;
			\item $\outcomeR$ otherwise.
		\end{itemize}
	\end{enumerate}
\end{proposition}

\begin{proof}
	We prove the statement for the case $a+b>c+d$, since the other case is symmetric. The case where $c=d=0$ is trivial; by convention, the first player cannot play, losing the game. Hence, let us assume that $c+d>0$.
	
	Note that, by Theorem~\ref{thm-moreHands}, Left wins if she is the first player. Furthermore, by Lemma~\ref{lem-lessHands}, Right can only win as the first player by removing a Left's hand with his move. Doing this would leave a position where both players have the same number of hands, and the determination of the outcome can be done by using Proposition~\ref{prop-magicNum2-sameNumberOfHands}.
	
	Now, the only way for Right to leave a $\outcomeP$-position is moving from an initial position $\handstwo{1^a,2}{1^a}$ with $a$ odd to $\handstwo{1^a}{1^a}$. Also, Right can leave an $\outcomeR$-position if he moves from $\handstwo{1^a,2^b}{1^{a+b-1}}$ with $b>1$ and $a$ odd to $\handstwo{1^a,2^{b-1}}{1^{a+b-1}}$. In all the other cases, playing first, Right leaves either an $\outcomeL$-position (if $b=0$, $d>0$ and $c$ is odd) or an $\outcomeN$-position (in all the other cases), and thus $P$ is an $\outcomeL$-position.
	
	Hence, for the position $P=\handstwo{1^a,2^b}{1^c,2^d}$, if $a+b>c+d$, then either $c=d=0$ and $P$ is a $\outcomeP$-position; or $d=0$, $b>0$ and $a$ is odd and $P$ is an $\outcomeN$-position; or, in all the other cases, $P$ is an $\outcomeL$-position. This finishes the proof.
\end{proof}

\subsection{Finger count 3} \label{subsec33}

A complete characterization of the outcomes of the \octopus/ positions with 3-fingered hands remains an open problem. Nonetheless, it is possible to determine the outcome of the initial position, a fact translated in the following theorem.

\begin{proposition}
	\label{prop-magicNum3}
	Let $a$ be a positive integer. The \octopus/ initial position $\hands{1^a}{1^a}{3}$ is an $\outcomeN$-position.
\end{proposition}

\begin{proof}
	Without loss of generality, we assume that Left is the first player. The proof is done by induction on the number of hands $a$. If $a=1$, then all the moves are forced, and Left wins:
	
	$$\hands{1}{1}{3}\xrightarrow{L}\hands{1}{2}{3}\xrightarrow{R}\hands{3}{2}{3}\xrightarrow{L}\{\,|\,\}.$$
	
	\vspace{0.3cm}		
	Assume now that $a>1$. From $P=\hands{1^a}{1^a}{3}$, Left can only move to the position $P'=\hands{1^a}{1^{a-1},2}{3}$. After that, $P'$ has two Right-options, and, so, we can consider the following two sequences of play:
	
	$$\hands{1^a}{1^a}{3}\xrightarrow{L}\hands{1^a}{1^{a-1},2}{3}\xrightarrow{R}\hands{1^{a-1},2}{1^{a-1},2}{3}\xrightarrow{L}\hands{1^{a-1},2}{1^{a-1}}{3}$$
	and
	$$\hands{1^a}{1^a}{3}\xrightarrow{L}\hands{1^a}{1^{a-1},2}{3}\xrightarrow{R}\hands{1^{a-1},3}{1^{a-1},2}{3}\xrightarrow{L}\hands{1^{a-1},3}{1^{a-1}}{3}.$$
	
	\vspace{0.3cm}
	In the first case, after Left's second move, Right is unable to eliminate any integer from Left's sequence, resulting in Left having more integers in her sequence when it is her turn. According to Theorem~\ref{thm-moreHands}, Left wins. In the second case, after Left's second move, there are two possible options for Right. The first one is to $\hands{1^{a-1}}{1^{a-1}}{3}$, which, by induction, is an $\outcomeN$-position. The second one is to $\hands{1^{a-2},2,3}{1^{a-1}}{3}$, where Left has more integers in their sequence than Right and have a winning move by Theorem~\ref{thm-moreHands}. In all possibilities, Left wins.
	
	We conclude that, playing first in $P$, Left wins. By symmetry, $P$ is a winning position for Right if he is the first player. Hence, $P$ is an $\outcomeN$-position.
\end{proof}

\section{Future work}\label{sec4}

A better understanding of the game values of $\hands{i}{j}{n}$ is clearly an interesting and challenging work to carry out. The computational complexity of \octopus/ (\emph{i.e.}, deciding which player has a winning strategy on a given position and with a given starting player) is also still open.

Also, as mentioned in the introduction, the \chopsticks/ variation that we studied is just the most obvious possibility. Indeed, here are several additional rules and variants that could be considered:
\begin{enumerate}
	\item The attack never changes; however there are two variants for how to remove a hand from play. The \emph{cut-off} variant, which we used, states that when an $n$-fingered hand should have more than $n$ fingers, then it is removed from play. The \emph{rollover} variant states that when an $n$-fingered hand with $x$ raised fingers is attacked by a hand with $y$ raised fingers, then it gets $x+y \bmod n$ raised fingers. So, for instance, with 5-fingered hands, attacking a hand with~4 raised fingers with a hand with~2 raised fingers will result in a hand with 1 raised finger. A hand is removed from play only if it reaches~0 raised fingers.
	\item The \emph{split} rule allows a player to shift fingers from one of their hands to another instead of attacking (note that a player may remove one of their own hands from play by doing so). However, simply reversing the two hands is not allowed.
	\item Expanding on this, the \emph{revive} rule allows a player to use a split move to revive one of their hands that had been removed from play. Variants of the revive rule may allow any revival, or restrict them to even splits (for instance, if a player has a hand with~4 raised fingers, then they can use it to revive another hand only by giving~2 fingers).
\end{enumerate}

\section{Acknowledgement}\label{sec5}

Antoine Dailly was supported by the International Research Center "Innovation Transportation and Production Systems" of the I-SITE CAP 20-25 and by the ANR project GRALMECO (ANR-21-CE48-0004).

Carlos Santos work is funded by national funds through the FCT - Funda\c{c}\~{a}o\linebreak para a Ci\^{e}ncia e a Tecnologia, I.P., under the scope of the projects\linebreak UIDB/00297/2020 and UIDP/00297/2020 (Center for Mathematics and \linebreak Applications).

\end{document}